\documentclass{article}
\usepackage[utf8]{inputenc}
\usepackage{bussproofs}
\usepackage{calc}
\usepackage{amsmath}
\usepackage{amsthm}
\usepackage{amssymb}
\usepackage{stmaryrd}
\usepackage{setspace}

\usepackage{MnSymbol}
\usepackage{cmll}
\usepackage{xcolor}

\newtheorem{definition}{Definition}

 \usepackage{url}
\usepackage[authoryear]{natbib}
\title{Game of Grounds}
\author{Davide Catta\\LIRMM - Montpellier University\\\texttt{davide.catta@lirmm.fr}\\ \\Antonio Piccolomini d'Aragona\\Aix-Marseille Univ, CNRS, Centre Gilles Gaston Granger\\Aix-en-Provence, France\\\texttt{antonio.piccolomini-d-aragona@univ-amu.fr}}
\date{}

\begin{document}

\maketitle

\begin{abstract}
\noindent In this paper, we propose to connect Prawitz's theory of grounds with Girard's Ludics. This connection is carried out on two levels. On a more philosophical one, we highlight some differences between Prawitz's and Girard's approaches, but we also argue that they share some basic ideas about proofs and deduction. On a more formal one, we sketch an indicative translation of Prawitz's theory grounds into Girard's Ludics relative to the implicational fragment of propositional intuitionistic logic. This may allow for a \emph{dialogical} reading of Prawitz's ground-theoretic approach. Moreover, it becomes possible to provide a formal definition of a notion of \emph{ground-candidate} introduced by Cozzo. 
\end{abstract}


\section{Introduction}

Dag Prawitz's \emph{theory of grounds} and Jean-Yves Girard's \emph{Ludics} are very recent semantics, shedding a new light upon some fundamental topics in contemporary mathematical logic.

Prawitz aims at explaining the compulsion exerted by proofs, in such a way that this epistemic power depends on the valid inferences of which a proof is built up. Girard, on the other hand, proposes a denotational semantics for Linear Logic, leading to a dialogical framework inspired by proof-search or game-theoretic approaches.

The two theories have different targets and employ different formal means. Furthermore, they structurally differ on many points, often fundamental. Nevertheless, they also share some common philosophical standpoints. First of all, the idea that evidence depends on (the possession of) objects constructed through primitive, meaning-constitutive or meaning-conferring, operations. Secondly, the idea that a proof is not an evidence object, but an act through which an evidence object is obtained. Thus, our first aim in this paper is that of providing a philosophical comparison between Prawitz's theory of grounds and Girard's Ludics, by highlighting their similarities, as well as their differences (Section 4).

In the light of these philosophical similarities, one could wonder whether the two frameworks can be linked also in a more formal sense. As a second point of our work, we suggest, in a purely indicative way, how the link can be found, by outlining a translation from Prawitz's theory of grounds into Girard's Ludics (Section 5). The link could allow for a dialogical reading of Prawitz's account. In addition, it could permit a rigorous description of the notion of ground-candidate introduced by Cesare Cozzo to fix some weak points of an earlier formulation of Prawitz's theory of grounds.

Before establishing the aforementioned philosophical and - indicative - formal links, however, we need to introduce the basic concepts of Prawitz's theory of grounds and of Girard's Ludics. Thus, this will be our starting point (Sections 2 and 3).

\section{Theory of grounds}

Through the theory of grounds (ToG), Prawitz aims at explaining how and why a correct deductive argument might compel us to accept its conclusion if we have accepted its premises/assumptions. Prawitz's previous semantics in terms of valid arguments and proofs (SAP)\footnote{See mainly \cite{Prawitz2008, Prawitz1974}.} had exactly the same purpose, but it suffered from some problems. Clearly, the reasons that led Prawitz to ToG are not the topic of this paper. However, some of them must be mentioned, because they will be important when discussing the relationship between Prawitz's and Girard's respective philosophical standpoints.

\subsection{From SAP to ToG}

Prawitz's first semantic proposal\footnote{See mainly \cite{Prawitz1974}.} is centered on the notion of valid argument. An argument is a sequence of first-order formulas arranged in tree form, where each node corresponds to an inference with arbitrary premises and conclusion. Validity of closed arguments - i.e. with no undischarged assumptions and unbound variables - is explained by saying that arguments ending with introductions, called canonical, are valid when their immediate sub-arguments are. Arguments ending with inferences in non-introductory form, called non-canonical, are instead valid when they can be reduced to a closed canonical valid argument. Open arguments are valid when all their closed instances are. Reductions employ constructive functions that transform trees ending with inferences in non-introductory form into trees where these inferences are removed. These functions are therefore a generalization of Prawitz's own reduction procedures for normalization in Gentzen's natural deduction systems.\footnote{See \cite{Gen34}, \cite{Pra65}.} In his later papers Prawitz proposes a similar approach, where the notion of valid argument in tree form is replaced by proofs understood as abstract objects built up of functions standing for valid inferences.\footnote{See mainly \cite{Prawitz2008}.} The definitions of notions such as closed canonical proof, closed non-canonical proof, and open proof, run basically like the corresponding ones in terms of valid arguments.

The question is now whether valid arguments and proofs, as defined in SAP, actually oblige us to accept their conclusion if we have accepted their assumptions. The epistemic power must stem from the valid inferences of which valid arguments and proofs are made up. Hence, what we need is an appropriate definition of the notion of valid inference. In SAP, this becomes the idea that an inference is valid if it gives rise to valid arguments and proofs when attached, respectively, to valid arguments and proofs. The overall framework, however, turns out to be unsatisfactory. One of the main problems stems from the interdependence of the concepts of valid inference and proof. The idea that proofs are chains of valid inferences seems to require a \emph{local} notion of valid inference; in SAP, instead, the order of explanation is reversed, for valid inferences are defined in the \emph{global} terms of the valid arguments and proofs they belong to. So, the compelling power cannot be explained by induction on the length of chains of valid inferences; since valid arguments and proofs may contain non-canonical steps, some valid inference may end a chain where a valid inference of the same kind, and of an equal or higher complexity may occur.\footnote{A similar point, although in a different context, is raised by \cite{tranchiniAphex} and \cite{Usberti2015}.}

Observe that the interdependence problem depends on the canonical/non-canonical distinction. So, a way out may be found by endorsing what is often called the \emph{proof-objects}/\emph{proof-acts} distinction introduced by Martin-L\"{o}f and developed by Sundholm\footnote{\cite{ML84,MartinLof1986ONTM}, \cite{sundholm}.}. The idea could be that proof-acts are chains of epistemic steps allowing to construct appropriate proof-objects. The former involve canonical or non-canonical moves, but we may require the latter to be always canonical, and specified by simple induction.\footnote{See also \cite{Tranchini2014b}.} Then, we may explain validity of inferences with respect to proof-objects, leaving proof-acts aside.

Clearly, valid inferences are exactly the steps that proof-acts are built up of. Hence, the above mentioned idea consists of explaining validity of inferences through the result that valid inferences yield - namely the proof-object they produce - rather than through how the result is obtained - namely the proof-act where they are used. This in turn means that what valid inferences generate is not only inferential structures, but abstract entities too. An inference cannot be the simple appending of a conclusion under some premises; also and above all, it has to be the application of an operation on proof-objects, with proof-objects as outputs in case of success.

As we shall see below, this strategy is exactly the one Prawitz follows with ToG.\footnote{A more detailed reconstruction of the content of this Section is in \cite{Piccolominithesis,piccolominidaragona2017}.}

\subsection{Grounds and their language}

The term ground is used by Prawitz to indicate “what a person needs to be in possession of in order that her judgment is to be justified or count as knowledge".\footnote{\cite{Prawitz2009}, 187.} According to this view, one should conceive of “evidence states as states where the subject is in possession of certain objects".\footnote{\cite{Prawitz2015}, 88.} But grounds also have to comply with a constructivist setup, whence they must be epistemic in nature. This means that they can be grasped according to the idea that “one finds something to be evident by performing a mental act".\footnote{\cite{Prawitz2015}, 88.} So, such an informal picture can be further specified by spelling out what kind of constructive operations build grounds. Hence, grounds will be what we may call operational entities.

A language of grounds refers to a background language. Endorsing the so-called \emph{formula-as-type conception}\footnote{See \cite{How69}. The \emph{formulas-as-types} conception is based on the idea that a formula should be understood as the class of its proofs, called its type. Thus, for example, a conjunction $A \wedge B$ corresponds to the cartesian product type $A \times B$, namely the class of pairs $\langle a, b \rangle$ with $a$ object in the type $A$ (proof of $A$) and $b$ object in type $B$ (proof of $B$); similarly, an implication $A \rightarrow B$ corresponds to the function space type $A \supset B$, namely, the class of functions $f(x^A)$ generating objects $f(g)$ in the type $B$ when applied to objects $g$ in the type $A$. The conception is consonant with the BHK interpretation of the meaning of the logical constants in terms of (canonical) proof-conditions of the formulas where such constants occur as main sign. In the case of an open formula $A(x_1, ..., x_n)$, the type is a function space, namely, the class of functions $f(x_1, ..., x_n)$ that produce objects in the type $A(k_1, ..., k_n)$ when applied to $n$ individuals. Since the output type now depends on the input type, we may speak with Martin-L\"{o}f of \emph{dependent types} (Martin-L\"{o}f 1985). One usually does the same with so-called hypothetical judgments $A_1, ..., A_n \vdash B$, where the input is given by objects in the types $A_1, ..., A_n$ and the output is an object in the type $B$.}, the formulas of the background language provide types for the terms - and other components - of the language of grounds, while the latter are meant to denote grounds for asserting the former - the assertion sign being the fregean $\vdash$. ToG is mainly concerned with first-order logic. So, hereafter, the background language will be a first-order language. As one usually does in intuitionistic frameworks, we put 

\begin{center}
$\neg A \stackrel{def}{=} A \rightarrow 0$
\end{center}
- where $0$ is an atomic constant for the absurd. Prawitz's main notions are relative to atomic bases. An \emph{atomic base} for a language of grounds over a given background language will be identified by individual and relational constants of the background language, plus a (possibly empty) set of atomic rules over the background language. Atomic rules are conceived of by Prawitz in the framework of so-called Post-systems, i.e. recursive sets of rules
\begin{prooftree}
\AxiomC{$A_1, ..., A_n$}
\UnaryInfC{$B$}
\end{prooftree}
over a first-order language $L$ such that: (1) $A_i$ and $B$ are atomic, and $A_i \neq 0$ ($1 \leq i \leq n$); (2) if $x$ occurs free in $B$, then there is $i \leq n$ such that $x$ occurs free in $A_i$.\footnote{Sometimes, one also considers atomic rules that bind individual variables or assumptions, but this point is not essential for our purposes.} Once a set of atomic rules has been set down, atomic derivations can be specified in a standard inductive way. 

We can outline a first example of language of grounds, labelled $\texttt{C}$, where terms are built up of symbols that correspond only to Gentzen's introduction rules. Given a first-order language $L$ and an atomic base $\mathfrak{B}$ over $L$ with atomic system $\texttt{S}$, the alphabet of $\texttt{C}$ contains names $c^A$ for atomic derivations of $A$ in $\texttt{S}$ with no undischarged assumptions and unbound variables, possibly indexed variables $\xi^A$ typed on formulas $A$ of $L$, typed operational symbols $k I$ with $k = \wedge, \vee, \rightarrow, \forall, \exists$, and a typed operational symbol $0_A$ for $A$ formula of $L$ expressing the explosion principle. Typed terms are defined in a standard inductive way, e.g. in the case of $\vee$, $\rightarrow$, $\exists$ and $0$ - we leave types of the operational symbols unspecified whenever possible -
\begin{itemize}
\item $T : A_i \in X \ \Rightarrow \ \vee I[A_i \vdash A_1 \vee A_2](T) : A_1 \vee A_2 \in X$ [$i = 1, 2$]
\item $T : B \in X \ \Rightarrow \ \rightarrow I \xi^A(T) : A \rightarrow B \in X$ [the typed-variable after the symbol indicates that this variable is bound by the symbol]
\item $T : A(t) \in X \ \Rightarrow \ \exists I[A(t) \vdash \exists x A(x)](T) : \exists x A(x) \in X$
\item $T : 0 \in X \ \Rightarrow \ 0_A(T) : A \in X$
\end{itemize}
With $\texttt{C}$ we can already state clauses fixing what counts as a ground for closed atomic formulas, for closed formulas with $\wedge$, $\vee$ and $\exists$ as main logical sign, and for $0$:

\begin{itemize}
\item[(A)] a ground over $\mathfrak{B}$ for atomic $\vdash A$ is any $c^A$
\item[($\wedge$)] a ground over $\mathfrak{B}$ for closed $\vdash A \wedge B$ is $\wedge I(T, U)$ where $T$ denotes a ground over $\mathfrak{B}$ for $\vdash A$ and $U$ denotes a ground over $\mathfrak{B}$ for $\vdash B$
\item[($\vee$)] a ground over $\mathfrak{B}$ for closed $\vdash A_1 \vee A_2$ is $\vee I[A_i \vdash A_1 \vee A_2](T)$ where $T$ denotes a ground over $\mathfrak{B}$ for $\vdash A_i$ with $i = 1$ or $i = 2$
\item[($\exists$)] a ground over $\mathfrak{B}$ for closed $\vdash \exists x A(x)$ is $\exists I [A(t) \vdash \exists x A(x)] (T)$ where $T$ denotes a ground over $\mathfrak{B}$ for $\vdash A(t)$ for some $t$
\item[($0$)] no $T$ denotes a ground over $\mathfrak{B}$ for $0$
\end{itemize}
On the contrary, the clauses for $\rightarrow$ and $\forall$ will involve total constructive functions of appropriate kind, taking individuals and/or grounds of given types as arguments, and producing grounds of given types as values. Like in the BHK-clauses\footnote{See for example Troelstra \& Van Dalen (1988). \nocite{Constructivism}}, the notion of total constructive function may be assumed as primitive:
\begin{itemize}
\item[($\rightarrow$)] a ground over $\mathfrak{B}$ for closed $\vdash A \rightarrow B$ is $\rightarrow I \xi^A(T)$ where $T$ denotes a ground over $\mathfrak{B}$ for $A \vdash B$, i.e. a constructive function over $\mathfrak{B}$ of type $A \vdash B$ [the typed variable after the symbol indicates that this variable is bound by the symbol]
\item[($\forall$)] a ground over $\mathfrak{B}$ for closed $\vdash \forall x A(x)$ is $\forall I x.(T)$ where $T$ denotes a ground over $\mathfrak{B}$ for $\vdash A(x)$, i.e. a constructive function over $\mathfrak{B}$ of type $A(x)$ [the individual variable after the symbol indicates that this variable is bound by the symbol]
\end{itemize}
The clauses can be considered as a ground-theoretic determination of the meaning of the logical constants. Therefore, the operational symbols of $\texttt{C}$ are primitive operations, and its terms may be qualified as canonical.

However, constructive functions have to be understood in an unrestricted way, and $\texttt{C}$ is too weak to express all of them. Some terms denote some constructive functions, but not all functions are denoted - e.g. there is no term denoting a constructive function of type $A_1 \wedge A_2 \vdash A_i$ ($i = 1, 2$). Thus, one must also bring in extensions of $\texttt{C}$ - in fact, as we shall see in Section 2.3, because of G\"{o}del's incompleteness no closed language of grounds permits to express all the grounds we need. The behaviour of an operation can be fixed through equations that show how to compute it on relevant values. As an example, consider the extension $\texttt{C}^*$ of $\texttt{C}$ whose alphabet contains new typed operational symbols $k E$ with $k = \wedge, \vee, \rightarrow, \forall, \exists$, standing for Gentzen's eliminations, and additional, inductively defined terms, e.g. in the case of $\vee$, $\rightarrow$ and $\exists$ - we leave types of the operational symbols unspecified -

\begin{itemize}
\item $T : A \vee B, U : C, V : C \in X \ \Rightarrow \ \vee E \ \xi^{A} \ \xi^{B}.(T, U, V) : C \in X$ [the typed-variables after the symbol indicates that these variables are bound by the symbol. Observe also that $U$ and $V$ are not necessarily of type $A \vdash C$ and $B \vdash C$ - these are just term-formation clauses. It is the equation below that permits us to interpret this operation semantically as one that yields a ground for $\vdash B$ when applied to grounds for $\vdash A \vee B$, $A \vdash C$ and $B \vdash C$ - think of the elimination rule for $\vee$ in Gentzen's natural deduction, and of the relative reduction as defined in Prawitz\footnote{See \cite{Pra65}.}]
\item $T : A \rightarrow B, U : A \in X \ \Rightarrow \ \rightarrow E(T, U) : B \in X$
\item $T : \exists x A(x), U : B \in X \ \Rightarrow \ \exists E \ x \ \xi^{A(x)}.(T, U) : B \in X$ [the typed- and individual variables after the symbol indicates that these variables are bound by the symbol. Observe also that $U$ is not necessarily of type $A(x) \vdash B$ - these are just term-formation clauses. It is the equation below that permits us to interpret semantically this operation as one that yields a ground for $\vdash B$ when applied to grounds for $\vdash \exists x A(x)$ and $A(x) \vdash B$ - think of the elimination rule for $\exists$ in Gentzen's natural deduction, and of the relative reduction as defined in Prawitz\footnote{See \cite{Pra65}.}]
\end{itemize}
The equations are standard conversions, e.g. in the case of $\vee E$, $\rightarrow E$ and $\exists E$

\begin{itemize}
\item $\vee E \ \xi^{A_1} \ \xi^{A_2}.(\vee I[A_i \vdash A_1 \vee A_2](T), U_1(\xi^{A_1}), U_2(\xi^{A_2})) = U_i(T)$
\item $\rightarrow E (\rightarrow I \xi^A (T(\xi^A)), U) = T(U)$
\item $\exists E \ x \ \xi^{A(x)}.(\exists I[A(t) \vdash \exists x A(x)](T), U(\xi^{A(x)})) = U(T)$
\end{itemize}

\noindent They show that the $kE$'s capture new total constructive functions, e.g. in the case of $\vee E$ and $\exists E$ of types

\begin{enumerate}
\item $A \vee B, (A \vdash C), (B \vdash C) \vdash C$
\item $\exists x A(x), (A(x) \vdash B) \vdash B$
\end{enumerate}
As a further example, consider the extension of $\texttt{C}^*$ obtained by adding a typed operational symbol $DS$ for disjunctive syllogism, with equations - we leave the type of $DS$ unspecified -

\begin{itemize}
\item $DS (\vee I[B \vdash A \vee B](T), U) = T$
\item $DS (\vee I[A \vdash A \vee B](T), U) = 0_B(\rightarrow E(T, U))$
\end{itemize}
At variance with $\texttt{C}$, the last two languages of grounds also contain non-canonical terms, that is, terms the outermost symbol of which is non-primitive.

Given a closed term $T$ in some language of grounds over $\mathfrak{B}$, we now say that $T$ denotes a ground over $\mathfrak{B}$ when it can be reduced to a term the outermost symbol of which is one of the operational symbols of $\texttt{C}$, that denotes a ground over $\mathfrak{B}$ according to one of the clauses ($\wedge$) - ($\forall$). Reduction employs the equations for the operational symbols of $T$, by replacing \emph{definiendum} by \emph{definiens}. For example
\begin{center}
$\rightarrow E (\rightarrow I \xi^{A \rightarrow A}_1(\vee \ \xi^{A \rightarrow A}_2 \ \xi^0(\vee I[A \rightarrow A \vdash A \rightarrow A \vee 0](\xi^{A \rightarrow A}_1), \xi^{A \rightarrow A}_2, 0_{A \rightarrow A}(\xi^0))), \rightarrow I \xi^A(\xi^A))$
\end{center}
by the equation for $\rightarrow E$ reduces to
\begin{center}
$\vee \ \xi^{A \rightarrow A}_2 \ \xi^0(\vee I[A \rightarrow A \vdash A \rightarrow A \vee 0](\rightarrow I \xi^A(\xi^A)), \xi^{A \rightarrow A}_2, 0_{A \rightarrow A}(\xi^0))$
\end{center}
which by the equation for $\vee E$ reduces to $\rightarrow I \xi^A(\xi^A)$. An open term denotes a constructive function over $\mathfrak{B}$ with domain $\mathfrak{D}$ and co-domain $\mathfrak{K}$ when all its closed instances denote a ground for $\mathfrak{K}$, where a closed instance is obtained by replacing individual variables with closed individual terms and typed variables with closed terms denoting grounds for the elements in $\mathfrak{D}$. So, if we replace $\rightarrow I \xi^A(\xi^A)$ with $\xi^{A \rightarrow A}_3$ in the example above, we obtain that the term we started from denotes a constructive function of type $A \rightarrow A \vdash A \rightarrow A$.\footnote{A better, but still partial development of the content of this Section is in \cite{piccolomini2018}; a more detailed development is in \cite{Piccolominithesis}.}

\subsection{Ground-theoretic validity}

Prawitz affirms that “to \emph{perform an inference} is, in addition to making an inferential transition, to apply an operation on the grounds that one considers oneself to have for the premises with the intention to get thereby a ground for the conclusion".\footnote{\cite{Prawitz2015}, 94.} An inference will be valid over $\mathfrak{B}$ when the application of the corresponding operation to the alleged grounds over $\mathfrak{B}$ for the premises actually yields a ground over $\mathfrak{B}$ for the conclusion. The inference is logically valid if it remains valid over all the bases. A proof (over $\mathfrak{B}$) can be now defined as a finite chain of valid (over $\mathfrak{B}$) inferences. ToG respects the proof-as-chains intuition: the notion of proof is non-circularly stated in terms of a \emph{local} notion of valid inference, in such a way that a proof involves only valid inferences. Crucially, “a proof of an assertion does not constitute a ground for the assertion but produces such a ground" \footnote{\cite{Prawitz2015}, 93.}, that is, ToG proofs are not objects but acts.

According to Prawitz's view, when an inferential agent carries a proof out, he/she comes into possession of a ground, not of a term. And grounds are, so to say, always canonical; terms describe how grounds are obtained - their being canonical or not depending on the kind of steps through which the possession is attained. Thus, when Gentzen's eliminations are linked to the non-primitive functional symbols of $\texttt{C}^*$, a performance of them on grounds for the premises will amount to nothing but a $\beta$-reduction. In more general cases, such as $DS$, we can in turn still resort to the old SAP idea of general procedures for obtaining canonical forms. Proofs can be conceived of as chains of applications of operations that, under relevant circumstances, coincide with computations on a sort of generalized non-normal forms.

\subsection{Cozzo's ground-candidates}

Cozzo moved four objections to the old formulation of ToG presented in Prawitz's first papers on the subject.\footnote{See \cite{Prawitz2009,Prawitz2012a}.} Only the fourth one is of interest for us here. Prawitz's definition of an act of inference was such that only valid inferences were inferences whereas, as Cozzo argues, “it seems reasonable to say that the experience of necessity of thought also characterizes the transition from mistaken premises devoid of corresponding grounds".\footnote{\cite{CozzoNecessity}, 114.} In replying to Cozzo's objections, Prawitz relaxes the definition of inference act by allowing what he calls \emph{alleged grounds}, i.e. “an [...] inference can err in two ways: the alleged grounds for the premisses may not be such grounds, or the operation may not produce a ground for the conclusion when applied to grounds for the premisses".\footnote{\cite{Prawitz2015}, 95.} However, Usberti observes that “Prawitz’s alleged grounds have nothing to do with ground-candidates: since no restriction is put on alleged grounds, they are entities of any kind. Now, an assertion based on an entity of any kind may be true or false, but it is difficult to see how it can be rational at all".\footnote{\cite{Usberti2017}, 525.} As Usberti observes, while Cozzo's first three objections might be overcome by Prawitz's adjustments, Cozzo's fourth objection remains unsolved.

As a remedy to his objections, Cozzo proposes an interesting notion of \emph{ground-candidate}, “a mathematical representation of the results of epistemic acts underlying mistaken premises". A ground-candidate “can be a genuine ground or a \emph{pseudo-ground}".\footnote{\cite{CozzoNecessity}, 114.} Can ground-candidates be characterized more precisely? As we shall see below, Girard's Ludics might suggest a promising account.

\section{Ludics}

Ludics was first proposed by Girard in his paper titled \emph{Locus Solum: from the rules of logic to the logic of rules}.\footnote{See \cite{Gir01locus}.} Its aim is that of studying the notions of proposition and proof (of type and element of a type) and of reconstructing them from a more primitive notion of interaction.

As is well known, the cut-rule allows for an ‘‘interaction" between two proofs. Given a proof of $A$ under hypotheses $\Gamma$, and a proof of $B$ under hypotheses $\Delta, A$, it yields a proof of $B$ under hypotheses $\Gamma, \Delta$. Gentzen's \emph{Haupsatz}  establishes that a derivation $\pi$ of $A$ under $\Gamma$ can be always ‘‘normalized" to a derivation $\pi'$ of $A$ under no more hypotheses than those in $\Gamma$, and where no cut-rule is employed. In the light of the Curry-Howard correspondence, normalization can be understood as the evaluation of a program to its normal form. 

Usually, the cut-rule is conceived as a deduction rule over a previously defined formal language. Once language and rules have been given, one can study the properties of cut-elimination. In order to attribute a computational meaning to proofs, a deterministic cut-elimination is required, i.e. given an instance of the cut-rule, there is only one way of eliminating it.

Ludics starts from an entirely opposite point of view. Given a deterministic cut-elimination defined on \emph{sui generis} computational objects, it seeks whether it is possible to recover formulas and rules. The objects of Ludics are an abstract counterpart of derivations in a well-behaved and complete fragment of multiplicative-additive Linear Logic. Such objects are built up of rules that ensure a deterministic cut-elimination on arbitrary numerical addresses. The latter represent the location that a formula may occupy in a derivation. This innovative approach, of which we shall say more below, is mainly inspired by a polarity phenomenon that we deal with in the next Section.

\subsection{Polarity}

Linear connectives are normally divided into two classes: of positive polarity - i.e. $\otimes$, multiplicative conjunction, $\oplus$, additive disjunction, $1$ multiplicative truth and $0$ additive false - and of negative polarity - i.e. $\with$, additive conjunction, $\parr$, multiplicative disjunction, $\top$, additive truth and $\bot$ multiplicative false. A formula is said to be negative (resp. positive) iff its main connective is negative (resp. positive). A linear connective $\star$ is said to be reversible iff, for each proof $\pi$ of a sequent $\Gamma, A$, for $A$ with main connective $\star$, there is a cut-free proof $\pi'$ of $\Gamma, A$, the last rule of which introduces $\star$.\footnote{See \cite{phdlaurent} for a detailed discussion of the phenomenon of polarity in Linear Logic.} The two notions are connected, in that negative connectives are reversible, whilst the positive ones are not.

Polarity and reversibility are very important for proof-search. Suppose we are looking for a proof of a sequent $\Gamma = \Gamma', A$, i.e. we are trying to prove $\Gamma$ by selecting a formula $A$ in $\Gamma$ and then by applying an inference on it. If $A$ is negative, the reversibility of its main connective ensures the existence of a cut-free proof of $\Gamma', A$, the last rule of which introduces $\star$. And this means that if, e.g., $A = B \parr C$, then we have only one bottom-up application of the rule, so that $\Gamma', B, C$ is provable too. Thus, selecting a negative formula in a sequent allows for an automatic proof-search procedure. But this strategy cannot be applied uniformly, because of the non-reversibility of positive connectives.

This notwithstanding, proof-search procedures can still be improved thanks to the following algorithm introduced by Andreoli.\footnote{See \cite{AndreoliFocusing}.} Given a sequent $\Gamma$: (1) if it contains negative formulas, focus on them and apply negative rules upwards until they are all decomposed, and there are no more negative formulas; (2) when you have finished, choose a positive formula randomly, and start decomposition by applying positive rules upwards, focusing at each step on its positive sub-formulas until you find negative formulas. Then, Andreoli proves that a sequent is derivable in Linear Logic, iff there is a focusing derivation of it, i.e. a derivation that can be built by applying his algorithm. It follows that one can restrict oneself to a set of derivations the elements of which consist of alternations of positive and negative steps only; positive steps cluster the positive branch of the algorithm - labelled (2) above - whereas negative steps cluster the negative branch of the algorithm - labelled (1) above. As an example of proof obtained through Andreoli's algorithm, consider the following one:

\begin{prooftree}
\AxiomC{}
\UnaryInfC{$\vdash A^\bot, A$}
\AxiomC{}
\UnaryInfC{$\vdash B^\bot, B$}
\BinaryInfC{$\vdash A , B, \mathbf{(A^\bot \otimes B^\bot)}$}
\UnaryInfC{$\vdash A , B, \mathbf{(A^\bot \otimes B^\bot)\oplus (A^\bot \otimes C^\bot)} $}
\AxiomC{}
\UnaryInfC{$\vdash A^\bot, A$}
\AxiomC{}
\UnaryInfC{$\vdash C^\bot, C$}
\BinaryInfC{$\vdash A , C, \mathbf{(A^\bot \otimes C^\bot)}$}
\UnaryInfC{$\vdash A , C,\mathbf{(A^\bot \otimes B^\bot) \oplus (A^\bot \otimes C^\bot)}$}
\BinaryInfC{$\vdash A , \mathbf{(B \with C)}, (A^\bot \otimes B^\bot) \oplus (A^\bot \otimes C^\bot)$}
\UnaryInfC{$\vdash \mathbf{A \parr (B\with C)} , \  (A^\bot \otimes B^\bot) \oplus (A^\bot \otimes C^\bot)$}
\end{prooftree}
Here, the bold formulas are those chosen for the step-wise proof-search procedure. The proof can be transformed into one where the two negative steps of the first sub-proof are merged into a single negative step, and where the same happens for the corresponding positive steps, i.e.
\begin{prooftree}
\AxiomC{}
\UnaryInfC{$\vdash A^\bot, A$}
\AxiomC{}
\UnaryInfC{$\vdash B^\bot, B$}
\BinaryInfC{$\vdash A, B, \mathbf{(A^\bot \otimes B^\bot) \oplus (A^\bot \otimes C^\bot)}$}
\AxiomC{}
\UnaryInfC{$\vdash A^\bot, A$}
\AxiomC{}
\UnaryInfC{$\vdash C^\bot, C$}
\BinaryInfC{$\vdash A , C, \mathbf{(A^\bot \otimes B^\bot) \oplus (A^\bot \otimes C^\bot)}$}
\BinaryInfC{$\vdash \mathbf{A \parr (B\with C)}, (A^\bot \otimes B^\bot) \oplus (A^\bot \otimes C^\bot)$}
\end{prooftree}
In this way, polarization permits to consider \emph{generalized connectives}, obtained by ‘‘merging" negative connectives with negative connectives, and positive connectives with positive connectives. In our example, the negative generalized connective is $- \parr (- \with -)$ and the positive generalized connective is $(- \otimes -) \oplus(- \otimes -)$. In addition, the sequents of the calculus end up having the form $\Gamma, \Delta$, with $\Gamma$ containing only positive formulas and $\Delta$ containing at most one negative formula. Although an arbitrary number of such connectives exists, all of them can be decomposed according to the same schemes of rules as follows - we use $P$ as a variable for positive formulas and $N$ as a variable for negative formulas:

\begin{prooftree}
\AxiomC{$\vdash N_{i_1},\Gamma_1 \quad \cdots \quad  \vdash  N_{i_{n_i}},\Gamma_n$}\RightLabel{positive rule}
\UnaryInfC{$\vdash(N_{1_1}\otimes \cdots \otimes N_{1_{n_1}})\oplus\cdots \oplus  (N_{p_1}\otimes\cdots \otimes  N_{p_{n_p}}),\Gamma_1,\ldots \Gamma_n $}
\end{prooftree}

\begin{prooftree}
\AxiomC{$\vdash P_{1_1}\cdots P_{1_n},\Gamma\quad\cdots \quad \vdash P_{1_{p_1}}\cdots P_{n_{p_n}},\Gamma $}\RightLabel{negative rule}
\UnaryInfC{$\vdash (P_{1_1}\parr\cdots \parr P_{1_n})\with\cdots \with (P_{1_{p_1}}\parr \cdots \parr P_{1_{p_n}}),\Gamma  
$}
\end{prooftree}
Thanks to the fact that De Morgan dualities hold for linear logic, and that the linear negation of a negative formula is a positive formula and vice versa, we can restrict our language to formulas composed by $\otimes, \oplus$ and write the schemes of rules above as follows:

\begin{prooftree}
\AxiomC{$P_{i_1}\vdash,\Gamma_1 \quad \cdots \quad  P_{i_{n_i}} \vdash \Gamma_n$}\RightLabel{positive rule}
\UnaryInfC{$\vdash(P^\bot_{1_1}\otimes \cdots \otimes P^\bot _{1_{n_1}})\oplus\cdots \oplus  (P^\bot_{p_1}\otimes\cdots \otimes  P^\bot_{p_{n_p}}),\Gamma_1,\ldots \Gamma_n $}
\end{prooftree}

\begin{prooftree}
\AxiomC{$\vdash P_{1_1}\cdots P_{1_n},\Gamma\quad\cdots \quad \vdash P_{1_{p_1}}\cdots P_{n_{p_n}},\Gamma $}\RightLabel{negative rule }
\UnaryInfC{$ (P^\bot_{1_1}\otimes\cdots \otimes P^\bot_{1_n})\oplus\cdots \oplus (P^\bot_{1_{p_1}}\otimes \cdots \otimes P^\bot_{1_{p_n}})\vdash \Gamma  
$}
\end{prooftree}


\subsection{From polarity to games}

Observe that each derivation in a calculus based upon the clustered rules of the previous Section is a tree in which positive rules are followed by negative rules and vice-versa - except the axiom rule. This allows us to consider ‘‘games" on a sequent.

Let us say that a \emph{move} is a couple $(F,\{F_1, ..., F_n\})$, where $F$ is a positive (negative) formula, called \emph{focus}, and $F_1, ..., F_n$ are some negative (positive) sub-formulas of $F$, called \emph{choices}. A move is of \emph{positive polarity} if its focus is a positive formula, and of \emph{negative polarity} otherwise. So, a \emph{game} is a non-empty sequence of moves such that: (1) two consecutive moves have opposite polarity; (2) the focus of a negative rule, except possibly the first, is one of the formulas in the choices of the positive move that immediately precedes it; (3) two distinct moves have distinct focuses. We finally define a strategy on $\Gamma \vdash \Delta$ to be a prefix-closed set of games on the same sequent. It is easily seen that the polarized clustered derivation in the previous Section can be considered as the set built up of the two games 

\begin{center}
    $\mathfrak{G}_1=(A \parr (B \with C),\{A,B\}),((A^\bot \otimes B^\bot)\oplus (A^\bot \otimes C^\bot),\{A^\bot,B^\bot\})$
\end{center}
\begin{center}
$\mathfrak{G}_2=(A \parr (B \with C),\{A,C\}),((A^\bot\otimes B^\bot)\oplus (A^\bot\otimes C^\bot),\{A^\bot,C^\bot\})$
\end{center}
The vice-versa holds too, i.e. the strategy can be converted into the proof above. 

Polarity is therefore useful to understand one of the main ideas that, as we shall see, inspire Girard's Ludics: proofs are strategies over a particular type of game. For this to make sense, though, we take a final step. The cut-rule

\begin{prooftree}
\AxiomC{$\vdash P$}
\AxiomC{$\vdash P^\bot$}
\BinaryInfC{$\vdash$}
\end{prooftree}
can be interpreted as indicating an interplay between a strategy $\pi$ for $P$ and an opposed strategy $\pi'$ for $P^\bot$. Clearly, if one between $P$ and $P^\bot$ is provable, the other cannot be so. Thus, if we want that every strategy ‘‘corresponds" to a proof in a polarized formal system as above, we need to admit proofs of both a proposition and its negation. To this end, we may introduce a new rule - call it \emph{daimon} - that permits to prove any possible sequent: 
 
 \begin{prooftree}
 \AxiomC{}\RightLabel{$\dagger$}
 \UnaryInfC{$\vdash \Gamma$}
 \end{prooftree}


\subsection{Ludics defined}

Ludics aims at overcoming the distinction between syntax (a formal system) and semantics (its interpretation). As is well-known, completeness implies that, for every $A$, either there is a proof of $A$ or a counter-model of $A$, i.e. a model of $\neg A$. This can be understood through a dialogical metaphor; in a two-persons debate, where one of the speakers tries to construct a proof of $A$ and the other tries to construct a counter-model of $A$, one and only one of the debaters can win. However, models and proofs are distinct entities and, in particular, there is no interaction between a proof of $A$ and a model of $\neg A$. Girard's idea is that of overcoming the syntax-semantic distinction by interpreting proofs with proofs, and counter-models with refutations. The properties of the attempted proofs of $A$ are verified by testing them, by means of a cut-elimination procedure, against attempted proofs of $\neg A$. Proofs should be understood as \emph{proof-search procedures} in the setup of the polarized system sketched in the two previous Sections. Given a certain conclusion, we try to guess the inference rule used to obtain it. If no rule is applied, we close the proof-search by ‘‘giving up", i.e. by a rule that encodes the information ‘‘I do not know how to keep on proving the conclusion, therefore I simply assert it".  This rule is what permits to have attempted proofs of both $A$ and $\neg A$. 


A Ludics derivation consists of rules extracted from those occurring in the above-told clustered derivations. Formulas are replaced by numerical addresses standing for the positions they may occupy during proof-search. In order to make this more precise, we now give a quick sketch of formalized Ludics.

We say that an \emph{address} is a string of natural numbers $i_1 \ \cdots \ i_n$. Two addresses are said to be \emph{disjoint} if none of them is a prefix of the other. A \emph{ramification} is a finite set of natural numbers. Given an address $\xi$ and a natural number $i$, with $\xi \ i$ we indicate the address obtained by putting $i$ at the end of $\xi$. Given an address $\xi$ and a ramification $I$, with $\xi \star I$ we indicate the set of the addresses $\{\xi \ i \ | \ \text{for every} \ i \in I\}$. A \emph{pitchfork} is an expression $\Gamma\vdash \Delta$ where $\Gamma,\Delta$ are finite sets of pairwise disjoint addresses such that $\Gamma$ contains at most one address. Addresses in $\Gamma$ are called \emph{negative}, while those in $\Delta$ are called \emph{positive}. A \emph{design} is a tree made of pitchforks, the last pitchfork being the \emph{base}, while the others are built through the \emph{rules}:\\

\noindent \textbf{daimon-rule} [already introduced in the previous Section]

\begin{prooftree}
\AxiomC{}
\RightLabel{$\dagger$}
\UnaryInfC{$\vdash \Gamma$}
\end{prooftree}

\noindent \textbf{Positive rule}

\noindent Let $I$ be a ramification and, for every $i \in I$, let the $\Gamma_i$ be pairwise disjoint and included in $\Gamma$. For every $i \in I$, the rule (finite) is

\begin{prooftree}
\AxiomC{$\cdots \ \xi \star i \vdash \Gamma_i \ \cdots$}
\RightLabel{($\vdash \xi, I$)}
\UnaryInfC{$\vdash \Gamma, \xi$}
\end{prooftree}

\noindent \textbf{Negative rule}

\noindent Let $N$ be a set of ramifications and, for every $I \in N$, let $\Gamma_I$ be included in $\Gamma$. For every $I \in N$, the rule (possibly infinite) is

\begin{prooftree}
\AxiomC{$\cdots \ \vdash \Gamma_I, \xi \star i \ \cdots$}
\RightLabel{($\xi \vdash N$)}
\UnaryInfC{$\xi \vdash \Gamma$}
\end{prooftree}

\noindent The design - not to be confused with the only rule applied in it -

\begin{prooftree}
\AxiomC{}
\RightLabel{$\dagger$}
\UnaryInfC{$\vdash \Gamma$}
\end{prooftree}
is called \emph{Daimon}. The only design not built by employing the rules above - called \emph{Fid} and responding to the idea of a positive conclusion with no rules above - is

\begin{prooftree}
\AxiomC{}
\RightLabel{$\Omega$}
\UnaryInfC{$\vdash \Gamma$}
\end{prooftree}

A \emph{cut} is given by an address that occurs once at the right of $\vdash$ (positive polarity) in the base of a design $\mathfrak{D}$, and once at the left of $\vdash $  (negative polarity) in the base of a design $\mathfrak{D}'$. Thus, a \emph{cut-net} is a finite non-empty set of designs such that: (1) the addresses in the bases are pairwise disjoint or equal; (2) every address occurs in at most two bases and, if it occurs in two bases, it occurs both positively and negatively; (3) the graph with vertices the bases and edges the cuts is connected and acyclic. The \emph{principal design} of a cut-net $\mathfrak{R}$ is the only design $\mathfrak{D} \in \mathfrak{R}$ with base $\Gamma \vdash \Delta$ such that $\Gamma$ is not a cut. The \emph{base} of a cut-net are the uncut addresses of the cut-net. 

We now proceed to an informal description of the process of normalization/interaction on \emph{closed} cut-nets, i.e. cut-nets the base of which is empty.\footnote{For a complete and formal definition, also on the open case, see \cite{Gir01locus}.} Given such a cut-net, the cut propagates over all the immediate sub-addresses as long as the action anchored on the positive pitchfork containing the cut corresponds to one of the actions anchored on the negative one. The process terminates when either the positive action anchored on the positive cut-fork is the daımon-rule, in which case we obtain a design with the same base as the starting cut-net, or no negative action corresponds to the positive one. In the latter case, the process is said to diverge. When normalization/interaction between two designs $\mathfrak{D}$ and $\mathfrak{D}'$ terminates and does not diverge, it ends up in Daimon, and $\mathfrak{D}$ and $\mathfrak{D}'$ are said to be orthogonal - indicated with $\mathfrak{D} \ \bot \ \mathfrak{D}'$. We give an example of terminating normalization on a  closed cut-net composed of two designs. Bold addresses are those through which the normalization/interaction procedures propagates.

\begin{prooftree}
\AxiomC{}
\RightLabel{$\dagger$}
\UnaryInfC{$\vdash \mathbf{\xi11}$}
\AxiomC{}\RightLabel{$\emptyset$}

\UnaryInfC{$\vdash \xi13$}
\BinaryInfC{$\mathbf{\xi1} \vdash$}
\UnaryInfC{$\vdash \mathbf{\xi}$}
\AxiomC{$\vdots$}
\noLine
\UnaryInfC{$\mathbf{\xi 11} \vdash$}
\UnaryInfC{$\vdash \mathbf{\xi1}$}
\AxiomC{$\vdots$}
\noLine
\UnaryInfC{$\vdash \xi2$}
\BinaryInfC{$\mathbf{\xi} \vdash$}
\noLine
\BinaryInfC{}
\end{prooftree}

\begin{prooftree}
\AxiomC{}
\RightLabel{$\dagger$}
\UnaryInfC{$\vdash \mathbf{\xi11}$}
\AxiomC{}\RightLabel{$\emptyset$}

\UnaryInfC{$\vdash \xi13$}
\BinaryInfC{$\mathbf{\xi1} \vdash$}
\AxiomC{$\vdots$}
\noLine
\UnaryInfC{$\mathbf{\xi11} \vdash$}
\UnaryInfC{$\vdash \mathbf{\xi1}$}
\noLine
\BinaryInfC{}
\end{prooftree}

\begin{prooftree}
\AxiomC{}
\RightLabel{$\dagger$}
\UnaryInfC{$\vdash \mathbf{\xi11}$}\AxiomC{$\vdots$}
\noLine
\UnaryInfC{$\mathbf{\xi11} \vdash$}
\noLine
\BinaryInfC{}
\end{prooftree}

\begin{prooftree}
\AxiomC{}
\RightLabel{$\dagger$}
\UnaryInfC{$\vdash$}
\end{prooftree}

\noindent Generalizing, we say that the \emph{ortoghonal} of a set $E$ of designs on the same base $\Gamma \vdash \Delta$ - written $E^\bot$ - is the set
\begin{center}
$\{\mathfrak{D} \ | \ \text{for all} \ \mathfrak{D}' \in E, \mathfrak{D} \ \bot \ \mathfrak{D}'\}$
\end{center}

For later purposes, we now need only four additional definitions. First of all, we say that a set of designs $E$ is a behaviour when $E = E^{\bot \bot}$. 
We give two simple examples of behaviours. Take any  design of the form  
\begin{prooftree}
\AxiomC{}
\RightLabel{$\emptyset$}
\UnaryInfC{$\vdash \xi $}
\end{prooftree}
 this type of design is called \emph{atomic bomb}, and take $E$ to be the singleton set containing the atomic bomb. The set $E^\bot$ only contains
\begin{prooftree}
\AxiomC{ }\RightLabel{$\dagger$}
\UnaryInfC{$ \vdash$}\RightLabel{$\emptyset$}
\UnaryInfC{$\xi\vdash$}
\end{prooftree}
This last design is also orthogonal to the Daimon
\begin{prooftree}
\AxiomC{}
\RightLabel{$\dagger$}
\UnaryInfC{$\vdash \xi$}
\end{prooftree}
Thus the set $E^{\bot\bot}$ will contain the atomic bomb and the Daimon - this behaviour represents the constant $1$. Consider now the following design, called \emph{skunk}. 
\begin{prooftree}
\AxiomC{}
\RightLabel{$\emptyset$}
\UnaryInfC{$\xi\vdash $}
\end{prooftree}
Its only orthogonal will be the the Daimon. Call $dai$ the set that only contains the Daimon. $dai^\bot$ will contain the skunk and any design of the form 
\begin{prooftree}
\AxiomC{$\vdots$}
\noLine
\UnaryInfC{}\RightLabel{$N$}
\UnaryInfC{$\xi\vdash$}
\end{prooftree}
where $N$ is a subset of the finite part of the power-set of the natural numbers. We denote this last behaviour by $\top$.

Before giving the second definition, observe that in the example of normalization above not all the addresses of the two designs that form the cut-net are explored by the normalization procedure, which means that in general, given a design $\mathfrak{D}$, only a sub-design $\mathfrak{D}'$ of $\mathfrak{D}$ is used to test a counter-design of opposite base. To get a clearer idea, consider the behaviour $\top$ defined above. The orthogonal of $\top$ is the behaviour $dai$. Given any design $\mathfrak{D}$ in $\top$, normalization between $\mathfrak{D}$ and the Daimon will not explore the addresses open by negative rules. So, the only design that is used in the interaction between an element of $\top$ and $\dagger$ is the skunk. This leads us to the following: given a design $\mathfrak{D}$ and a behaviour $G$, we call \emph{incarnation} of $\mathfrak{D}$ in $G$ - indicated with $|\mathfrak{D}|_G$ -  the smallest sub-design $\mathfrak{D}'$ of $\mathfrak{D}$ which is still in $G$.
As a third definition, we say that $\mathfrak{D}$ is \emph{material} in a behaviour $G$ when $\mathfrak{D} = |\mathfrak{D}|_G$. The skunk in $\top$, as well as the two designs in $1$ are material. Finally, given a behaviour $G$, we call \emph{incarnation} of $G$ - indicated with $|G|$ - the set
\begin{center}
    $\{\mathfrak{D} \ | \ \mathfrak{D} \in G \ \text{and} \ \mathfrak{D} = |\mathfrak{D}|_G\}$
\end{center}

To conclude, two points must be underlined. First of all, $\dagger$ stands for paralogism; it captures the idea of abandoning the dialogue or game on a position that one is not able to justify further. When interaction yields $\dagger$, one of the two designs is a locally winning strategy; though, the same strategy may lose in other contexts. A globally winning strategy is instead one that can never be defeated; a design can thus be understood as a proof when the interaction with any of its ortoghonals gives $\dagger$. It follows that designs are not necessarily proofs - they may be nothing but attempted proofs. Secondly, designs are cut-free; cuts only occur when interaction is defined. The idea is that a type is a set of cut-free attempted proofs; more precisely, a type amounts to a set of cut-free paraproofs that behave in the same way in interactions with orthogonal designs. More precisely, in Ludics a type is nothing but a behaviour, as the latter notion has been defined above.

\section{Differences and similarities}

We now propose a philosophical comparison between ToG and Ludics. Although some deep differences may be detected, we shall argue that these two theories share equally deep tenets, somehow inspiring the general framework they provide.

\subsection{Differences: order, types and bidirectionalism}

The first and most striking difference between ToG and Ludics concerns the logic which they aim at interpreting. Girard's Ludics can be considered as an interpretation of second-order multiplicative-additive Linear Logic with weakening. Every formula $A$ is interpreted as a behaviour, and this permits to prove that: (1) if $\pi$ is a proof of $A$, then there exists a $\dagger$-free, material design $\mathfrak{D}$ in the behaviour associated to $A$ such that $\mathfrak{D}$ is the interpretation of $\pi$; (2) if $\mathfrak{D}$ is a material and $\dagger$-free design in a behaviour $A$, then $\mathfrak{D}$ is the interpretation of a proof $\pi$ of $A$. Prawitz's ToG, instead, mainly aims at interpreting first-order intuitionistic logic. As is well-known, a Linear interpretation of first-order intuitionistic logic requires the modality operator $!$, allowing for contraction. So, for example, the intuitionistic $A \rightarrow B$ can be defined as $(!A) \multimap B$. In the version of Ludics we have been referring so far, no interpretation of the type $!A$ is given. However, such interpretation \emph{can} be given, for example in \emph{$C$-Ludics}\footnote{See \cite{Terui}.}, or in \emph{Ludics with repetitions}.\footnote{See Faggian \& Balsadella (2011) \nocite{faggian-basaldella}.} A difference between $C$-Ludics and Ludics with repetitions is that the designs of $C$-Ludics may contain cuts, whereas those in Ludics with repetitions - as already in Girard's original formulation - are always cut-free. Since the cut-free character of designs will be one of the common points we will highlight between ToG and Ludics, the comparison we propose should be thought as referring to Ludics with repetitions, rather than to $C$-Ludics.

This first kind of difference also concerns the order on which Prawitz and Girard reason. Prawitz's grounds and terms of a language of grounds are thought of as referring to a \emph{first-order} background language. The latter provides types for either classifying abstract objects, or labelling syntactical expressions that denote such objects. Higher-order theories of grounds may be of interest, but they would face the same difficulties as those met in any other constructivist approach. In fact, it is well-known that $n$th-order logics for $n \geq 2$ imply a loss of what Dummett \cite{DummettSeas} called molecularity on introduction rules.\footnote{See \cite{DummettSeas}, \cite{CozzoTeoria}.} So, a ground for a second-order existential $\vdash \exists X \alpha(X)$ should be defined through a primitive operation, say $\exists^2 I$, by requiring that $\exists^2 I(\mathfrak{U}, T, X)$ denotes a ground for $\vdash \exists X \alpha(X)$ iff $T$ denotes a ground for $\vdash \alpha(\mathfrak{U})$; but $\mathfrak{U}$ may contain $\exists X \alpha(X)$ as a sub-formula, and the definition could be no longer compositional. In addition, higher-order constructivist setups might suffer from impredicativity and paradoxical phenomena.\footnote{See \cite{PISTONE_2018}.} On the contrary, Ludics is a \emph{second-order} theory, and this would appear as an unsurmontable obstacle between Prawitz's and Girard's approaches.

Recent works, however, have shown that Ludics is suitable for first-order quantification\footnote{See Fleury \& Quatrini (2004) \nocite{firstorderLudics}.}, or for first-order Martin-L\"{o}f dependent types\footnote{See \cite{Sironi,sironithesis}.}. Moreover, the propositional level is much less problematic - although, as we will see, not entirely unproblematic. The difference concerning logics is therefore undoubtedly important, but not as worrying as it may appear at a first glance.

The second difference is what we may call the \emph{typed vs untyped} distinction. Inspired by the Curry-Howard isomorphism, ToG adopts the so-called \emph{formulas-as-types} conception, in the light of which, as already remarked, a background language provides \emph{types} for grounds and terms. Conversely, Ludics is fully \emph{untyped}. Nonetheless, Ludics aims at recovering types out of a more primitive notion of interaction. In fact, we \emph{do have} types herein; they are reconstructed as sets of designs respecting certain constraints, namely, sets of designs equal to their bi-orthogonal, namely, behaviours. By attributing a pivotal role to the notion of interaction, Ludics focuses on typability, rather than on typedness; metaphorically, it seeks the conditions under which a dialogue obtains, as well as the conditions under which a dialogue can be said to employ logical means. So, while in ToG types are primitive sets of proofs, in Ludics they are non-primitive sets of well-behaving designs. And since a design can be understood as a paraproof, proofs will only constitute the \emph{subset} of a type - i.e. the set of the designs of a behaviour that win in every possible interaction.

The third and deeper difference between Prawitz's and Girard's approaches is finally the following one; while ToG is clearly a verificationist semantics, in that it bases the explanation of meaning on primitive operations that mirror Gentzen's introduction rules, Ludics seems to be more akin to bidirectionalism.

\emph{Verificationism} can be roughly described as the idea of explaining meaning in terms of the \emph{conditions} for asserting correctly; its dual is \emph{pragmatism}, according to which meaning should be fixed in terms of the correct \emph{consequences} of an assertion. Gentzen's natural deduction provides a paradigmatic picture for both verificationism and pragmatism; the former is centred on introduction rules, with respect to which elimination rules have to be justified, whereas the latter proceeds the other way round.

\emph{Bidirectionalism} is instead the standpoint according to which meaning must be explained in terms of \emph{two} primitive notions; the conditions for asserting correctly, and the conditions under which one can safely modify the assumption of an assertion.\footnote{For the Linear Logic framework, see \cite{Zeilberger}; for the natural deduction framework, see \cite{SchroederHeister2009SequentCA}; an interesting discussion of bidirectionalism in connection with grounding can be found in \cite{Francez2015pts} - Francez's grounding is however entirely different from Prawitz's one.} Not surprisingly, bidirectionalism fits better with Gentzen's sequent calculus that comes with two kinds of introduction rules - right introductions and left introductions.

As remarked by Schroeder-Heister, left introductions for the logical constant $k$ can be understood as generalized elimination rules where a major premise with main sign $k$ occurs as an assumption, and replaces assumptions of lower complexity, which may in turn be discharged. Now, it can be shown that negative ‘‘clustered" rules, as discussed in Section 3.1, fit with generalized elimination rules of this kind. To verify why Ludics is more akin to bidirectionalism, it is sufficient to recall that its rules abstract from the ‘‘clustered" ones; \emph{positive} Ludics rules mirror right introductions, whilst \empty{negative} Ludics rules mirror left introductions, i.e. generalized eliminations. Crucially, both positive and negative rules are \emph{primitive}, relative to the determination of types.

\subsection{Similarities: objects/acts and computation}

Despite the more or less relevant differences, it seems to us that ToG and Ludics share at least two basic ideas. These tenets somehow inspire, from a philosophical as well as from a formal standpoint, the overall framework which Prawitz and Girard envisage for a rigorous reconstruction of deduction. Both such ideas can be understood through the lens of the aforementioned distinction between \emph{proof-objects} and \emph{proof-acts}.

As previously written, the distinction is explicitly at play in ToG; on the one hand, we have grounds - what one is in possession of when justified - and, on the other, we have proofs - acts by means of which grounds are obtained. These two sides, which we may call respectively \emph{abstract} and \emph{operational}, interact via the terms of the languages of grounds; a term denotes a ground, and codes a proof delivering the denoted ground.

Also Ludics involves an abstract and an operational side. Here, the objects are attempted proofs - designs - while the dynamics of deduction is represented by interaction - cut-elimination in cut-nets. Undoubtedly, Girard proposes a very peculiar standpoint about the act of (para)proving, i.e. a dialogical or game-theoretic one. Nevertheless, there is still the idea that a (para)proof-act produces a (para)proof-object; interaction looks for the (locally or globally) winning design and, if converging, it results in a normal form representing a strategy that the (locally or globally) winner may endorse for being (locally or globally) justified.

Now, both in ToG and in Ludics, the abstract/operational articulation relies upon what seems to be the same programmatic idea - being the first similarity we are going to point out. Prawitz's objects are, as we have seen, always \emph{canonical}, since grounds are always obtained by applying primitive operations; ToG acts, instead, may be \emph{canonical or not}, depending on how the ground they yield is obtained. Likewise, Girard's objects, i.e. designs, are always \emph{cut-free}; \emph{cuts} only occur in cut-nets, i.e. in interaction.

Thus, to Prawitz's abstract level, inhabited by canonical objects, corresponds Girard's abstract level, inhabited by normal designs. Indeed, canonicity can be looked at as a semantic generalization of normal form, where derivations are taken not as elements generated by a fixed formal system, but as broad structures through which one defines notions such as validity and consequence.\footnote{See \cite{sch2006} about the so-called \emph{fundamental corollary} of Prawitz's normalization theory.} To Prawitz's operational level, where we do have a distinction between canonical and non-canonical cases, corresponds Girard's operational level, where we do have cuts. Once again, we have a semantic generalization of cut-elimination; to justify a non-canonical step, one requires harmony with respect to the canonical cases, and this is done by showing how maximal peaks that the non-canonical step would create may be appropriately dropped out.

Reasonably, one could now wonder where such a symmetry stems from. To answer this question, we have to turn to the second similarity between ToG and Ludics. Recall that Prawitz's notion of ground accounts for evidence. But justification takes root in meaning. Hence, grounds are built up of \emph{primitive} and \emph{meaning-constitutive} operations. Within Ludics, the parallel idea is that (tentative) evidence, represented by designs, depends on some \emph{primitive} rules, the interaction of which in cuts allows for \emph{determination of meaning} of definable types. Non-canonicity and cuts appear as soon as we focus, not on what evidence is, but on how evidence is acquired. This may be because for both Prawitz and Girard a proof is an act of a very special kind: it is a \emph{computation} that removes non-primitive elements.

To be more precise, we know that Prawitz takes inferences as applications of operations on grounds. Given the way operations on grounds are understood, i.e. as fixed through equations that show how to compute non-canonical terms, we have pointed out that, under relevant circumstances, inferences could be understood as \emph{generalized reduction steps}. Analogously, in Girard's framework, interaction is a process through which a (locally or globally) winning design may be found at the end of a dialogue or game between opposite positions. The act of proving may be therefore understood as exactly this proof-search, proceeding via \emph{cut-elimination}.

However, the parallelism between the abstract sides of Prawitz's and Girard's theories undergoes some restrictions. In fact, the claim that to Prawitz's objects - the grounds - correspond Girard's object - the designs - must be accompanied by the observation that designs are not objects on the same level or of the same nature as grounds. A ground reifies an evidence state for a \emph{specific} judgement or assertion; as such, it is of a \emph{specific} type, according to the proposition or sentence involved in the judgment or assertion for which it is a ground. On the other hand, a design reifies the \emph{moves} in a proof-game, or in proof-search, independently from judgements or assertions at issue in the interplay; as already remarked, Ludics is untyped precisely because it aims at reconstructing the pure \emph{dynamics} of giving and asking for reasons, and at recovering typed logical strategies out of a more basic notion of interaction. To use a slogan, we could say that, while Prawitz's objects are objects in the ‘‘Fregean" sense of being saturated entities obtained by filling unsaturated entities, Girard's objects are ‘‘interactional" objects.

On the other hand, Girard's designs come conceptually closer to Prawitz's grounds when they are looked at within the global context of a behaviour, i.e. when they share stable interactional properties with other designs. Once a class of strategies has been, so to say, closed under cut-elimination, we are allowed to speak about a type. Typing is thus obtained by abstracting from overall properties that objects of the same type are expected to show when used in (attempted) deduction. This is also the linchpin of our indicative proposal for formally linking Prawitz's ToG and Girard's Ludics.

\section{Grounds in Ludics}

After the philosophical comparison between ToG and Ludics, we seek whether these two theories can also be compared on a more formal level. We remark, however, that what follows is just an introductory framework, requiring further work and refinement.

\subsection{A translation proposal: the implicational fragment}

Before starting, we pay our due to Sironi.\footnote{See \cite{Sironi,sironithesis}.} When proposing an understanding of Prawitz's grounds in Girard's framework, we will largely rely upon her embedding of Martin-L\"{o}f's type theory in Ludics.

By limiting ourselves to the implicational fragment of intuitionistic logic, we show in a purely indicative way how Prawitz's grounds for formulas of the kind $A \rightarrow B$ could be translated onto Girard's framework as designs respecting certain constraints - we discuss only in the conclusive remarks to what extent the translation may apply to the other intuitionistic first-order constants.

The leading idea of our indicative mapping is that of undertaking Girard's own line of thought, and looking at a type $A$ as a behaviour $G^A$. A ground for $\vdash A$ is to be a $\dagger$-free element of $|G^A|$, the incarnation of $G^A$ . Finally, a cut-net $\mathfrak{R}$ between two appropriate $\dagger$-free elements of incarnations, and such that $[[\mathfrak{R}]]$ is a $\dagger$-free element of $|G^A|$, will here stand for a non-canonical term denoting the ground that corresponds to $[[\mathfrak{R}]]$. We now want to add something to justify the chosen strategy.\footnote{Sironi's setup is more complex. Let us sketch it quickly. Given an address $\xi$, a ramification $I$ and a set of ramifications $N$, a \emph{positive action} is either $(+, \xi, I)$ or $\dagger$; a \emph{negative action} is $(-, \xi, N)$. A \emph{chronicle} $\mathfrak{c}$ is a non-empty, finite, alternate sequence of actions such that: (1) each action of $\mathfrak{c}$ is initial or justified by a previous action of opposite polarity; (2) the actions of $\mathfrak{c}$ have distinct addresses; (3) if present, $\dagger$ is the last action of the chronicle. A $\dagger$-\emph{shorten} of a chronicle $\mathfrak{c}$ is $\mathfrak{c}$ or a prefix of $\mathfrak{c}$ ended by $\dagger$. The $\dagger$-\emph{shortening} of $E$ - written $E^{\dagger}$ - is the set of designs obtained from $E$ by $\dagger$-shortening chronicles. $E$ is \emph{principal} iff its elements are $\dagger$-free and $|E^{\bot \bot}| = E^{\dagger}$. Sironi takes: a type $A$ to be $(P^A)^{\bot \bot}$, with $P^A$ principal; a canonical element of a type to be an element of $P^A$; a cut-net $\mathfrak{R}$ such that $[[\mathfrak{R}]] \in P^A$ to be a non-canonical element denoting the canonical element that $[[\mathfrak{R}]]$ corresponds to. The appeal to principal set of designs is due to the fact that Sironi aims at having types generated by their canonical elements; indeed, $P^A$ is a kind of minimal generator of $(P^A)^{\bot \bot}$.}

First of all, why should one understand a type as a behaviour? A behaviour is a set of designs equal to its bi-orthogonal. This morally means that a behaviour contains all the necessary and sufficient information for a dialogue or game to take place. Thus, a behaviour yields meaningfulness, in the sense that it permits to punctually counter-argue in any possible way over its base.

From this point of view, a ground is a $\dagger$-free element of the incarnation of a behaviour for three reasons. It is $\dagger$-free because, of course, one expects a ground not to involve any paralogism. On the other hand, the incarnation of a behaviour is the subset of the behaviour the elements of which actually operate in interactions. Hence, although we could have demanded a ground to belong to the behaviour as such, we more specifically require a ground to belong to its incarnation so to have something minimal that is used in dialogues or games. Finally, a $\dagger$-free element $\mathfrak{D}$ of the incarnation $|G|$ of a behaviour $G$ (if any) enjoys the following property: given $\mathfrak{D}' \in |G^\bot|$, $[[\mathfrak{D}, \mathfrak{D}']] =$ Daimon, and since $\mathfrak{D}$ is $\dagger$-free, $\dagger$ must occur in $\mathfrak{D}'$. Hence, we have what we required of designs-as-proofs in Section 3.3; $\mathfrak{D}$ wins in every possible interaction with elements of the incarnation of the orthogonal of the behaviour which it belongs to. Observe that, if $\mathfrak{D}' \in G^\bot - |G^\bot|$, $[[\mathfrak{D}, \mathfrak{D}']]$ may diverge; this simply means that $\mathfrak{D}'$ is, so to say, not answering the question raised by $\mathfrak{D}$, i.e. we are not in the presence of an actual dialogue or game.

In the specific case of implication, we thus proceed as follows. We know that a ToG ground over $\mathfrak{B}$ for $\vdash A \rightarrow B$ is $\rightarrow I \xi^A(T)$, where $T$ denotes a constructive function over $\mathfrak{B}$ of type $A \vdash B$. Suppose that we have appropriately translated $\mathfrak{B}$ onto the Ludics framework.\footnote{An example of how atomic types can be put in Ludics is given again by \cite{Sironi,sironithesis} for the types $\mathbb{N}$ and $\texttt{List}$.} Suppose also that the types $A$ and $B$ have been inductively determined as behaviours $G^A$ and $G^B$. Observe finally that, as already remarked, interaction, i.e. normalization of cut-nets, is a deterministic procedure.

Let us indicate with $|G^A|_F$ and $|G^B|_F$ the sets of the $\dagger$-free elements of $|G^A|$ and $|G^B|$ respectively. Given $\mathfrak{D} \in |G^A|_F$, and given $\mathfrak{D}'$ such that $[[\mathfrak{D}, \mathfrak{D}']] \in |G^B|_F$, we can take $[[\mathfrak{D}, \mathfrak{D}']]$ as the result of a constructive function of type $A \vdash B$. Hence, we put

\begin{center}
$A \rightarrow B = \{\mathfrak{D}' \ | \ \text{for every} \ \mathfrak{D} \in |G^A|_F, \ [[\mathfrak{D}, \mathfrak{D}']] \in |G^B|_F\}^{\bot\bot}$.
\end{center}
We must take the bi-orthogonal, because we have no guarantee that the above mentioned set is a behaviour. On the other hand, it holds that, for every set of designs $E$, $E^\bot$ is a behaviour. Such technical move has no disturbing effect. We can understand the bi-orthogonal as the interaction-closure of the set as such. Furthermore, the elements of $A \rightarrow B$ are designs of base $\alpha \vdash \beta$, where $\vdash \alpha$ is the base of $G^A$ and $\vdash \beta$ is the base of $G^B$.\footnote{Strictly speaking, what we have defined here is the linear arrow $\multimap$, which is usually done by putting $A \multimap B = A^\bot \parr B$. However, our definition is equivalent to the standard one, except that it takes the advantage of not passing through $\parr$ - which would have required the introduction of many new notions and definitions. Observe also that the idea of defining a type as a set closed under bi-orthogonality is quite standard in Linear Logic frameworks, in particular in the Geometry of Interaction program, as remarked in Naibo, Petrolo \& Seiller (2016) \nocite{naibo:halshs-01313400}.}

So, $\rightarrow I \xi^A(T)$ can be understood as a $\dagger$-free element $\mathfrak{D}'$ of $|A \rightarrow B|$. Observe that, for every $\mathfrak{D} \in |G^A|_F$, inductively corresponding to a ground $g$ over $\mathfrak{B}$ for $\vdash A$, the cut-net with $\mathfrak{D}$ and $\mathfrak{D}'$ can be taken as the non-canonical term $\rightarrow E (\rightarrow I \xi^A(T), g)$, denoting a ground over $\mathfrak{B}$ for $\vdash B$. The application of $\mathfrak{D}'$ to $\mathfrak{D}$ produces an interaction $[[\mathfrak{D}, \mathfrak{D}']]$ that ends in an element of $|G^B|_F$.

We quickly outline a concrete - although very simple - example of translation. Consider the orthogonal of the behaviour $\top$ introduced in section 3.3, i.e. the behaviour containing only the Daimon. Call it $0$. Consider now the following recursive design (called $Fax$): 

\begin{prooftree}
\AxiomC{$ \ \ \ \ \ \quad \quad \vdots \ Fax_{\xi.i\vdash \xi.i}$}
\noLine
\UnaryInfC{$\cdots \qquad \xi'.i\vdash \xi.i \qquad \cdots $}\RightLabel{$(\xi',I)$}
\UnaryInfC{$ \cdots \qquad \vdash  \xi.I,\xi' \qquad \cdots $}\RightLabel{$(\xi,\mathcal{P}_f(\mathbb{N}))$}
\UnaryInfC{$ \xi\vdash \xi'$}
\end{prooftree}
$Fax$ is  the Ludics interpretation of the identity axiom in sequent calculus. It is a $\dagger$-free  element of the incarnation of any behaviour $A \rightarrow A$. In particular, it is an element of the behaviour $0\to 0$. In fact  given the only material design in $0$ with base $\xi$ - observe that this design \emph{is not} $\dagger$-free - the normalization between $Fax$ and this latter gives a Daimon based on $\xi'$, which is again a material design in $0$. We thus interpret the term $\rightarrow I \xi^0 (\xi^0)$ of Tog as $Fax$ in $(0 \rightarrow 0)$. A similar translation may be applied - with some more machinery - to any $A$, taking a $\dagger$-free material element of the incarnation of $A$ as inductively corresponding to a ground for $\vdash A$.

However, we remark that, because of the first difference we have highlighted in Section 4.1 between ToG and Ludics, the indicative mapping we have proposed can work only for \emph{linear} terms of ToG, i.e. terms where every occurrence $\rightarrow I$ binds \emph{exactly one} typed-variable. An extension should take into account what we have already said in the same Section about intuitionistic implication. That is, we need a Ludics interpretation of the modality operator $!$, and then we can take $A \rightarrow B$ as $(! A) \multimap B$. Since we want designs to be cut-free, this should be done following the work of Faggian and Basaldella\footnote{See Faggian \& Basaldella (2011). \nocite{faggian-basaldella}} rather then Terui with $C$-Ludics.\footnote{See \cite{Terui}.}

If one accepts our reconstruction above of a ToG type as a behaviour, and of a ToG ground as a $\dagger$-free element of the incarnation of a behaviour, one also easily notices that a type contains more than grounds. Grounds are designs that win in every relevant interaction, but the behaviour they belong to may contain designs with paralogisms or designs losing in some interactions. Thus, we could simply define Cozzo's ground-candidates, discussed in Section 2.5, as generic elements of a behaviour which are pseudo-ground if they are not $\dagger$-free or do not belong to the incarnation of the behaviour. In this way, a ground-candidate would be a structure representing one's strategy in a dialogue or game, which might win in some contexts, but lose in others.

\subsection{Cozzo's ground-candidates reconsidered}

Let us now turn back to Cozzo's fourth objection to the old formulation of ToG discussed in Section 2.4. Prawitz's 
definition implied that only valid inferences were inferences, and according to Cozzo this is misleading in that necessity of thought, what Prawitz aims at explaining, can be experienced also in inferences with mistaken premises. To take into account this phenomenon, Cozzo suggests the introduction of ground-candidates, representing a reasoning with possibly mistaken premises, so that a ground-candidate is either an actual ground, or just a pseudo-ground. Cozzo's view seems to involve three aspects that we had better distinguish now:
\begin{enumerate}
    \item an inference can be performed on wrong premises;
    \item an inference performed on wrong premises can be valid;
    \item a valid inference on wrong premises can be epistemically compelling.
\end{enumerate}
Point 1 requires that inferential operations are defined, not only on grounds as Prawitz's operations illustrated in Section 2.2, but also on pseudo-grounds, and hence, more in general, on ground-candidates. Thus, as in point 2, the inference can be valid even if what one is in possession of are just pseudo-grounds. What matters is that, \emph{when applied to grounds} for the premises, the operation always produces a ground for the conclusion.

However, for such an inference to be also compelling, as point 3 requires, it must hold that the inferential operation still produces ground-candidates for the conclusion when applied to pseudo-grounds for the premises. The result of the application might be a ground, but we cannot expect it to be so in every case. And if it were not at least a pseudo-ground, i.e. an epistemic support based on which one feels entitled - perhaps wrongly - to assert the conclusion, there would be no way to maintain that the inference has a compelling character. The question is now whether ground-candidates, in particular pseudo-grounds, can be adequately characterised within the framework of ToG. We shall argue that they cannot.

The problem with ToG is that entities that \emph{are not} grounds cannot reasonably be considered also as pseudo-grounds. This is particularly evident in the case of ground-terms $T$ that denote something which \emph{is not} a ground for $A \vdash B$ - call it $f(\xi^A)$. All we know is that, for some ground $g$ for $\vdash A$, $f(g)$ is not a ground for $\vdash B$. However, for this to be sufficient to conclude that $f(\xi^A)$ is at least a pseudo-ground for $A \vdash B$, it seems plausible to require that $f(g)$ is at least a pseudo-ground for $\vdash B$.

Now, a possibility that may not be excluded is simply that $f(g)$ either diverges, i.e.
\begin{center}
    $f(g) = f_1(f(g)) = f_2(f_1(f(g))) = ... = f_n(...(f_2(f_1(f(g)))...) = ...$
\end{center}
or gives rise to a loop, i.e.
\begin{center}
    $f(g) = f_1(f(g)) = f_2(f_1(f(g))) = ... = f_1(f_2( ... ((...(f_2(f_1(f(g))))...)) = f(g)$.
\end{center}
For example, it has been suggested\footnote{See \cite{Tennant1995,tennant}. For a critical discussion, see Petrolo \& Pistone (2018), \nocite{petroloPistone} \cite{Tranchini2018}.} that this is what happens with paradoxes in a proof-theoretic framework, but without going that far we may consider a ‘‘ping-pong" where $f(\xi^A)$ is defined by the equation
\begin{center}
    $f(g) = f_1(g, g)$
\end{center}
and $f_1(\xi^A, \xi^A)$ by the equation
\begin{center}
     $f_1(g, g) = f(g)$.
\end{center}
It is clear that we are not entitled to look at $f(g)$ as a pseudo-ground for $\vdash B$. In fact, the computation yields no value, and a computation that yields no value cannot constitute an epistemic support for - possibly mistaken - assertions. As a consequence, a \emph{modus ponens} from $A \rightarrow B$ and $A$ to $B$ where one is in possession of $\rightarrow I \xi^A(T)$ and $g$, represented by $\rightarrow E (\rightarrow I \xi^A(T), g)$, would produce - according to the definition of $\rightarrow E$ - the divergent or looping computation $f(g)$. So, in spite of its validity, it would not be compelling. We may postulate, of course, that $f(g)$ is an alleged ground for $\vdash B$, as Prawitz does in response to Cozzo's objection.\ But, as Usberti remarks, this alleged ground is an entity of any kind, and an assertion based on an entity of any kind can be hardly understood as rational.

Observe that, in the example of divergent or looping computation provided above, we are not only unable to hit upon a ground for $\vdash B$. We are even unable to obtain a \emph{canonical} object of type $B$. Thus, an obvious way out may be that of requiring that, for it to be a pseudo-ground for $A \vdash B$, $f(\xi^A)$ yields in all cases canonical objects of type $B$, although for some ground for $\vdash A$ the canonical object produced is not a ground, but only a pseudo-ground for $\vdash B$ - where (some or all) the arguments to which the primitive operation is applied are pseudo-grounds. This should be generalized to ground-candidates, by requiring that a ground-candidate for $A \vdash B$ is a constructive function that, for every ground-candidate for $\vdash A$, yields a canonical object of type $\vdash B$. The notion of ground-candidate can be specified by induction on the complexity of formulas, in the same way as the notion of ground. For example, in the case of implication the clause would run as follows:
\begin{itemize}
    \item[($\rightarrow^*$)] $T$ denotes a ground-candidate for $A \vdash B$ iff $\rightarrow I \xi^A(T)$ denotes a ground-candidate for $\vdash A \rightarrow B$.
\end{itemize}
Accordingly, we should modify our notion of inference, by requiring that operations on grounds applied in inference steps are only defined on ground-candidates of the kind just defined. The \emph{modus ponens} we have considered above, thus, should not be accepted as an inference, since in that case $\rightarrow I$ has an immediate sub-argument that stands for a divergent or looping computation.

In this way, we have excluded from the class of ground-candidates all those entities that are not grounds, and that involve divergent or looping computations. ToG might be modified so as to fit with this idea, but it is doubtful whether this can be done without deeply modifying Prawitz's original project. However, this may constitute the topic of further works. For the moment, we remark that, morally, the result we aim to is exactly what we obtain by understanding Prawitz's grounds as $\dagger$-free elements of the incarnation of a behaviour and Cozzo's ground-candidates as simple elements of a behaviour. If a design $\mathfrak{D}$ belongs to a behaviour $A$, it encodes enough information for it to be tested against any design in the behaviour $ A^\bot$. In particular, this means that $\mathfrak{D}$ has the ‘‘shape" of a normal derivation of $A$, independently of whether such derivation actually exists. The inferential pattern that it is built up is made of primitive steps, independently of whether such steps can be understood as drawn from initial sequents in a specific formal calculus. Finally, it is part of its being the element of a behaviour that any interaction of $\mathfrak{D}$ with any element of $A^\bot$ does not give rise to divergent or looping computations. The interaction \emph{produces} a value which may be a ground, but which is not necessarily so. Even in this case, however, the value is still an element of a behaviour, whence it is typed.

\section{Conclusion}

The philosophical and (indicative) formal bridging between Prawitz's theory of grounds and Girard's Ludics we have proposed seems to allow for a dialogical reading of the former. There is still the idea of explaining evidence in terms of some primitive operations that yield normal/canonical objects. Furthermore, evidence is thought of as obtained by performing constructive acts which correspond to reduction/cut-elimination over non-primitive steps. However, evidence now stems from interaction, and since none of the interacting agents might be right, the game may end up in a ground-candidate that happens not to be an actual ground. A picture of this kind seems to be difficult to be accomplished within ToG. As we have seen, attributing types to pseudo-grounds in ToG may be a very difficult task, which requires deep modifications in Prawitz's original project. A type for a pseudo-ground would be nothing but a formal label, and this may not be enough to ensure that the object is constructed by using inferences that are meaning-constitutive, or at least meaning-justifiable. This situation is \emph{a priori} excluded in Ludics, where we can define a pseudo-ground for $A$ as an element of the behaviour representing $A$. Elements of the behaviour are abstract counterparts of cut-free derivations for $A$, and this independently of whether \emph{any} such derivation exists.

Admittedly, our first step in relating ToG and Ludics is very limited. In particular, in order to relate the two theories in a more credible way, one should be able to regain the full "power" of intuitionistic operations into Ludics. We are quite confident that this can be done in the frame proposed by Faggian and Basaldella\footnote{See Faggian \& Basaldella (2011). \nocite{faggian-basaldella}}. Even if Ludics with repetitions does not enjoy all the properties of the version of Ludics we have referred to here, it still does enjoy the ones that are central for our philosophical discussion.

\medskip

\bibliographystyle{chicago}
\bibliography{biblio}

\end{document}